\documentclass{article}   	
\usepackage{amsmath,amsthm,amssymb,bm,graphicx,caption,float,lscape}   
\usepackage[hyphens]{url}
\usepackage[colorlinks,citecolor=blue,urlcolor=blue]{hyperref}
\pdfoutput=1
\def\P{{\rm P}}        
\def\E{{\rm E}}        
\def\Var{{\rm Var}}    
\def\SD{{\rm SD}}      

\theoremstyle{plain}
\newtheorem{theorem}{Theorem}
\newtheorem{lemma}[theorem]{Lemma}

\theoremstyle{remark}
\newtheorem*{remark}{Remark} 

\allowdisplaybreaks

\newenvironment{newreferences}
               {\section*{References}
                \raggedright
                \begin{list}{}{\setlength{\itemsep}{0pt}
                               \setlength{\parsep}{0pt}
                               \setlength{\labelwidth}{0pt}
                               \setlength{\leftmargin}{12pt}
                               \setlength{\labelsep}{0pt}}
                \setlength{\itemindent}{-12pt}
               }{\end{list}}

\title{Does the first-serving team have \\ a structural advantage in pickleball?}
\author{Daryl R. DeFord\thanks{\,Department of Mathematics and Statistics, Washington State University. Email: \href{mailto:daryl.deford@wsu.edu}{daryl.deford@wsu.edu}.}\ \ and Stewart N. Ethier\thanks{\,Department of Mathematics, University of Utah. Email: \href{mailto:ethier@math.utah.edu}{ethier@math.utah.edu}.}}
\date{}

\begin{document}

\maketitle

\begin{abstract}
In pickleball doubles with conventional side-out scoring, points are scored only by the serving team. The serve alternates during a game, with each team serving until it has faulted twice, except at the beginning of the game, in which case the first-serving team serves until it has faulted once. A game to $n$ can be modeled by a Markov chain in a state space with $4n^2+10$ states. Typically, $n=11$ or $n=15$. The authors, both pickleball players, were motivated by the question in the title. Surprisingly, the answer to that question depends on the number of points needed to win. In a game to 11, the first-serving team has a very slight disadvantage, whereas, in a game to 15, the first-serving team has a very slight advantage. It should be noted that these advantages and disadvantages are so small that they cannot be detected by simulation and are revealed only by an analytical solution.  The practical implication is that a team that is offered the choice of side or serve should probably choose side. 

We investigate the probability of winning a game to 11, as well as the mean and standard deviation of the duration (or the number of rallies) of a game to 11.  We compare these results with the corresponding ones when modified rally scoring is used in a game to 21.  We also investigate the title question for a hybrid form of rally scoring that combines modified rally scoring and traditional doubles server rotation.\medskip

\noindent\textit{Keywords}: pickleball, side-out scoring, modified rally scoring, Markov chain, fundamental matrix, absorption probabilities, computer algebra
\end{abstract}

\section{Introduction}

Pickleball is a tennis-like sport that has recently seen a surge in popularity.  In 2023 it was named the fastest-growing sport in the United States for the third consecutive year by the Sports \& Fitness Industry Association.  A recent estimate put the number of US participants as high as 36 million (Golden, 2023).  Like tennis it can be played in singles or doubles format, but doubles is more widely played, so we restrict our attention to that format.  Additionally, the scoring system of  pickleball singles is equivalent to an obsolete volleyball scoring system, which was studied by Calhoun, Dargahi-Noubary, and Shi (2002), while pickleball doubles presents new complications that have not been fully analyzed.

The scoring system in pickleball is called \textit{side-out scoring}.  Points are scored only by the serving team.  The serve alternates during a game, with each team serving until it has faulted twice, except at the beginning of the game, in which case the first-serving team serves until it has faulted once.  The server (or referee if there is one) calls the score before each serve: $i$-$j$-$k$ with $i$ being the serving team's score, $j$ being the receiving team's score, and $k$ being 1 or 2 to to indicate whether it is the serving team's first or second server.  The game begins at 0-0-2.

Let $p_A$  (resp., $p_B$) be the probability that Team $A$ (resp., Team $B$) wins a rally when serving.  Assuming rallies are independent, we can calculate such quantities as the probability that Team $A$ wins a game to 11 when Team $A$ is first server, as a function of $p_A$ and $p_B$.  The so-called two-bounce rule in pickleball prohibits the serving team from volleying the return of serve, and this allows the receiving team to go on offense first, giving them a small advantage in each rally.  We expect that, at least in high-level play, $p_A$ and $p_B$ lie in the interval $[0.4,0.5]$.  An estimate based on data from 40 pro gold medal matches was given by Xu (2022) as $p_A=p_B\approx0.44$.  (In particular, the interesting portion of the parameter space in pickleball is very different from that in tennis, where the server has a large advantage, and both are very different from that in volleyball, where the serving team has a large disadvantage.)

The problem we want to solve is related to the analogous problem in tennis, whose solution is well known (Kemeny and Snell, 1976, \S 7.2).  But there are several reasons why the pickleball problem is more complicated.

\begin{enumerate}
\item \textit{Game to 11 or 15}.   A tennis game is played to 4.  A pickleball game is played to 11 or 15.  In either game, one must win by two or more.

\item \textit{Alternating serves}.  In a tennis game, the same player (or team) serves throughout.  In pickleball, the serve alternates during the game.  

\item \textit{Side-out scoring}.  In pickleball, points are scored only by the serving team, a system known as side-out scoring.  Tennis uses rally scoring.

\item \textit{Two faults before side out}.  In pickleball doubles, a team serves until it has faulted twice (or just once to begin the game).  This is in contrast to pickleball singles.
\end{enumerate}

\begin{figure}
\centering
\includegraphics[width=3.5in]{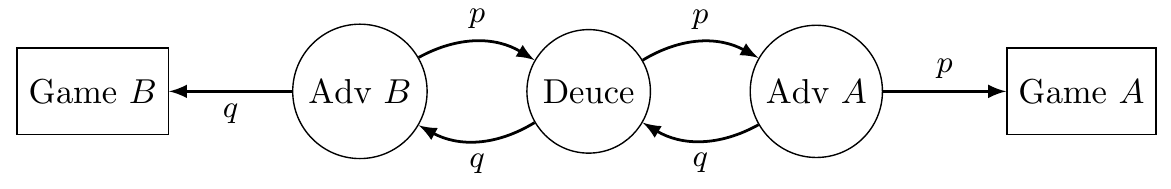}
\caption{\label{fig:tennis-deuce}The deuce game at tennis.}
\end{figure}

To emphasize these distinctions, consider a tennis game from deuce.  With $p$ being the probability that the server wins a point, the game is mathematically analogous to a gambler's ruin Markov chain on $\{0,1,2,3,4\}$.  See Figure~\ref{fig:tennis-deuce}.

\begin{figure}
\centering
\includegraphics[width=4.75in]{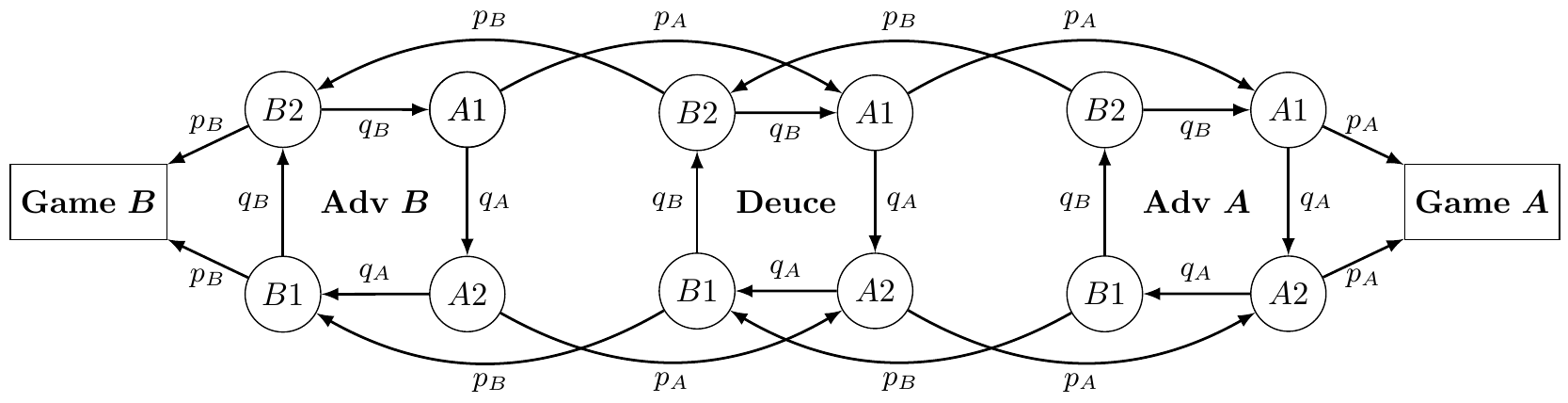}
\caption{\label{fig:pickleball-deuce}The deuce game at pickleball.}
\end{figure}

We can similarly consider a pickleball game from deuce.\footnote{Pickleball does not use this term, but it could.}  Here we need two parameters, $p_A$, the probability that Team $A$ wins a point when serving, and $p_B$, the probability that Team $B$ wins a point when serving.  But now the state space has 14 states, not just 5.  See Figure~\ref{fig:pickleball-deuce}.

Assuming rallies are independent, we can model a full game to $n$ by a Markov chain in a state space with $4n^2+10$ states. Typically, $n=11$ or $n=15$, so the one-step transition matrix is large but sparse (at most two nonzero entries in each row).  Application of standard textbook results allows us to compute the probability that Team $A$ wins the game, either from 0-0-2 or from any score, as well as the mean and variance of the game's duration (as measured by the number of rallies), again either from 0-0-2 or from any score.

The title question can be answered by evaluating the advantage of the first-serving team, namely
\begin{align*}
&\P(\text{\textnormal{Team $A$ wins game to} $n$}\mid \text{\textnormal{Team $A$ serves first}}) \\
&\quad{}-\P(\text{\textnormal{Team $A$ wins game to} $n$}\mid \text{\textnormal{Team $B$ serves first}}),
\end{align*}
at least for $n=11$ and $n=15$.  A somewhat analogous quantity in tennis, namely
\begin{align*}
&\P(\text{\textnormal{Player $A$ wins set}}\mid \text{\textnormal{Player $A$ serves first}}) \\
&\quad{}-\P(\text{\textnormal{Player $A$ wins set}}\mid \text{\textnormal{Player $B$ serves first}}),
\end{align*}
is known to be 0, that is, a tennis set is \textit{service neutral} (Newton and Pollard, 2004; MacPhee, Rougier, and Pollard, 2004; Newton and Keller, 2005).  This again assumes independence of rallies. 

In pickleball we find that, for the values of $p_A$ and $p_B$ that would be typical of high-level play, the first-serving team has a very slight disadvantage when $n=11$ and a very slight advantage when $n=15$.  But these advantages and disadvantages are so small that they cannot be detected by simulation and are revealed only by an analytical solution.  The practical implication, as first pointed out by Xu (2022), is that a team that is offered the choice of side or serve should probably choose side.  Although the side advantage (due primarily to wind and sun) can be mitigated by a side change when one of the teams reaches $(n+1)/2$ points, there may be a nontrivial advantage in having the better side last.

There has been discussion in the pickleball community of whether pickleball should switch from side-out scoring to rally scoring, as volleyball did in 1999.  Independent of that, a form of scoring called \textit{modified rally scoring} has been adopted by Major League Pickleball, a pro league created in 2021.  In rally scoring, a point is awarded to the winner of each rally, and that team then serves to begin the next rally.  In a game to $n$, modified rally scoring coincides with rally scoring until one of the teams has reached $n-1$ points.  From there, that team can score additional points only when serving.  Furthermore, given that one team has reached $n-1$ points, the opposing team, once it reaches $n-3$ points, can score additional points only when serving.  It turns out that, in modified rally scoring to 21, the first-serving team has an advantage if $p_A>q_B:=1-p_B$ and a disadvantage if $p_A<q_B$, as might be expected.

Here we compare results for modified rally scoring to 21 with those for side-out scoring to 11.  We also investigate the title question for a hybrid form of rally scoring that combines modified rally scoring and traditional doubles server rotation.

\section{Methodology}

By modeling a game as a Markov chain, we can readily compute such quantities as win probabilities and the mean and variance of game duration.  We begin with a one-step transition matrix $\bm P$ for a Markov chain in the state space $\{1,2,\ldots,m+2\}$.  States $1,2,\ldots,m$ are assumed transient and states $m+1$ and $m+2$ are assumed absorbing.  The one-step transition matrix $\bm P$ can be expressed in block form as 
\begin{equation}\label{block-P}
\bm P=\begin{pmatrix}\bm Q & \bm S\\ \bm 0 & \bm I_2\end{pmatrix},
\end{equation}
where $\bm Q$ is the $m\times m$ submatrix of $\bm P$ corresponding to one-step transitions from transient states to transient states, $\bm S$ is the $m\times 2$ submatrix of $\bm P$ corresponding to one-step transitions from transient states to absorbing states, $\bm 0$ is a $2\times m$ matrix of 0s, and $\bm I_2$ is the $2 \times 2$ identity matrix.  The assumption that a state $i\in\{1,2,\ldots,m\}$ is transient means that, from state $i$, state $m+1$ or $m+2$ can be reached in a finite number of steps, that is, $(\bm P^n)_{i,m+1}>0$ or $(\bm P^n)_{i,m+2}>0$ for some $n\ge1$.

We denote the Markov chain with one-step transition matrix $\bm P$ by $\{X_n,\; n=0,1,\ldots\}$, and we denote the time to absorption by
\begin{equation*}
T:=\min\{n\ge0: X_n=m+1\text{ or }X_n=m+2\}.
\end{equation*}
Let $\bm I$ denote the $m\times m$ identity matrix, and denote by $\bm1$ the column vector of 1s of length $m$.  The quantities of interest are given in the following lemma.

\begin{lemma}[Kemeny and Snell, 1976, Theorems 3.3.7 and 3.3.5]\label{methodology}
Under the above assumptions on $\bm P$, the matrix $\bm I-\bm Q$ is nonsingular.  Let $\bm M:=(\bm I-\bm Q)^{-1}$ denote the fundamental matrix.  Given a transient state $i\in\{1,2,\ldots,m\}$,
\begin{equation*}
\P_i(X_n=m+1\text{\textnormal{ for all $n$ sufficiently large}})=(\bm M\bm S)_{i,m+1},
\end{equation*}
\begin{equation*}
\E_i[T]=(\bm M\bm 1)_i, \quad\text{and}\quad \Var_i(T)=[(2\bm M-\bm I)\bm M\bm 1]_i-[(\bm M\bm 1)_i]^2.
\end{equation*}
\end{lemma}

The subscript $i$ on $\P$, $\E$, and $\Var$ indicates the initial state.  The lemma generalizes easily to the case in which the Markov chain has not just an initial state $i$ but an initial distribution $\bm\pi$.

\section{Side-out scoring}

In pickleball doubles with conventional side-out scoring, points are scored only by the serving team.  The serve alternates during a game, with each team serving until it has faulted twice, except at the beginning of the game, in which case the first-serving team serves until it has faulted once.  A game to $n$ can be modeled by a Markov chain in a state space with $4n^2+10$ states.  $4n^2+8$ of the states are of the form $(i,j,k)$, where $(i,j)\in(\{0,1,2,\ldots,n-1\}\times\{0,1,2,\ldots,n-1\})\cup\{(n-1,n),(n,n-1)\}$ and $k\in\{1,2,3,4\}$.  In state $(i,j,k)$, Team $A$ has $i$ points, Team $B$ has $j$ points, and $k$ indicates whether the server is Team $A$'s first server (1), Team $A$'s second server (2), Team $B$'s first server (3), or Team $B$'s second server (4).  There are also two absorbing states, ``Team $A$ wins'' and ``Team $B$ wins''.  The three scores, $(n-1,n-1)$, $(n-1,n)$, and $(n,n-1)$, can be thought of, in the language of tennis, as ``deuce'', ``advantage Team $B$'', and ``advantage Team $A$''.\footnote{Pickleball does not use these terms.  This is merely a device for ensuring a finite state space, and it involves no loss of generality.  There is a further possible simplification, corresponding to regarding 30-30 as deuce in a tennis game.  We choose not to adopt this simplification because it makes the labeling of states less straightforward.}   The transition matrix $\bm P$ can be written in block form as in \eqref{block-P}, where $m=4n^2+8$.

Let $p_A$ (resp., $p_B$) be the probability that Team $A$ (resp., Team $B$) wins a rally when serving, let $q_A:=1-p_A$ (resp., $q_B:=1-p_B$), and define the entries of $\bm Q$ and $\bm S$ as follows.

For $i,j\in\{0,1,\ldots,n-2\}$,
\begin{small}
\begin{align*}
Q((i,j,1),(i+1,j,1))&=p_A, & Q((i,j,1),(i,j,2))&=q_A, \\
Q((i,j,2),(i+1,j,2))&=p_A, & Q((i,j,2),(i,j,3))&=q_A, \\
Q((i,j,3),(i,j+1,3))&=p_B, & Q((i,j,3),(i,j,4))&=q_B, \\
Q((i,j,4),(i,j+1,4))&=p_B, & Q((i,j,4),(i,j,1))&=q_B.
\end{align*}
\end{small}
For $i\in\{0,1,\ldots,n-2\}$,
\begin{small}
\begin{align*}
Q((i,n-1,1),(i+1,n-1,1))&=p_A, & Q((i,n-1,1),(i,n-1,2))&=q_A, \\
Q((i,n-1,2),(i+1,n-1,2))&=p_A, & Q((i,n-1,2),(i,n-1,3))&=q_A, \\
S((i,n-1,3),\text{Team }B\text{ wins})&=p_B, & Q((i,n-1,3),(i,n-1,4))&=q_B, \\
S((i,n-1,4),\text{Team }B\text{ wins})&=p_B, & Q((i,n-1,4),(i,n-1,1))&=q_B.
\end{align*}
\end{small}
For $j\in\{0,1,\ldots,n-2\}$,
\begin{small}
\begin{align*}
S((n-1,j,1),\text{Team }A\text{ wins})&=p_A, & Q((n-1,j,1),(n-1,j,2))&=q_A, \\
S((n-1,j,2),\text{Team }A\text{ wins})&=p_A, & Q((n-1,j,2),(n-1,j,3))&=q_A, \\
Q((n-1,j,3),(n-1,j+1,3))&=p_B, & Q((n-1,j,3),(n-1,j,4))&=q_B, \\
Q((n-1,j,4),(n-1,j+1,4))&=p_B, & Q((n-1,j,4),(n-1,j,1))&=q_B.
\end{align*}
\end{small}
Finally,
\begin{small}
\begin{align*}
Q((n-1,n-1,1),(n,n-1,1))&=p_A, & Q((n-1,n-1,1),(n-1,n-1,2))&=q_A, \\
Q((n-1,n-1,2),(n,n-1,2))&=p_A, & Q((n-1,n-1,2),(n-1,n-1,3))&=q_A, \\
Q((n-1,n-1,3),(n-1,n,3))&=p_B, & Q((n-1,n-1,3),(n-1,n-1,4))&=q_B, \\
Q((n-1,n-1,4),(n-1,n,4))&=p_B, & Q((n-1,n-1,4),(n-1,n-1,1))&=q_B, \\
\noalign{\medskip}
Q((n-1,n,1),(n-1,n-1,1))&=p_A, & Q((n-1,n,1),(n-1,n,2))&=q_A, \\
Q((n-1,n,2),(n-1,n-1,2))&=p_A, & Q((n-1,n,2),(n-1,n,3))&=q_A, \\
S((n-1,n,3),\text{Team }B\text{ wins})&=p_B, & Q((n-1,n,3),(n-1,n,4))&=q_B, \\
S((n-1,n,4),\text{Team }B\text{ wins})&=p_B, & Q((n-1,n,4),(n-1,n,1))&=q_B, \\ 
\noalign{\medskip}
S((n,n-1,1),\text{Team }A\text{ wins})&=p_A, & Q((n,n-1,1),(n,n-1,2))&=q_A, \\
S((n,n-1,2),\text{Team }A\text{ wins})&=p_A, & Q((n,n-1,2),(n,n-1,3))&=q_A, \\
Q((n,n-1,3),(n-1,n-1,3))&=p_B, & Q((n,n-1,3),(n,n-1,4))&=q_B, \\
Q((n,n-1,4),(n-1,n-1,4))&=p_B, & Q((n,n-1,4),(n,n-1,1))&=q_B.
\end{align*}
\end{small}

All entries not specified are 0.  For specificity, it is convenient to order the states $(i,j,k)$ lexicographically, followed by states ``Team $A$ wins'' and ``Team $B$ wins''.  That completes the specification of the transition matrix $\bm P$.  The initial state is $(0,0,2)$ if Team $A$ is first server, $(0,0,4)$ if Team $B$ is first server.

\begin{theorem}\label{Thm:windiff-so}
Consider a game of pickleball doubles with conventional side-out scoring.  Let $p_A$ (resp., $p_B$) be the probability that Team $A$ (resp., Team $B$) wins a rally when serving, and assume that rallies are independent and $p_A+p_B>0$. Then
\begin{align*}
f_n(p_A,p_B)&:=\P(\text{\textnormal{Team $A$ wins game to} $n$}\mid \text{\textnormal{Team $A$ serves first}})\\
&\qquad{}-\P(\text{\textnormal{Team $A$ wins game to} $n$}\mid \text{\textnormal{Team $B$ serves first}})
\end{align*}
is a symmetric rational function on $[0,1]\times[0,1]-\{(0,0)\}$.  

In the special case in which the teams are evenly matched, that is, $p_A=p_B=:x\in(0,1]$,
\begin{equation*}
f_n(x,x)=2\,\P(\text{\textnormal{first-serving team wins game to} $n$})-1.
\end{equation*}
\end{theorem}

\begin{remark}
The function $f_n$ might be called the \textit{first-serving team's advantage} under conventional side-out scoring to $n$.  A negative advantage is of course a disadvantage.  Formulas for $f_{11}$ and $f_{15}$ are deferred to the Appendix.  Xu (2022), using iteration in a spreadsheet, found the approximate numerical value of $f_{11}(0.4,0.4)$, accurate to about five decimal places.
\end{remark}

\begin{proof}
With $p_A+p_B>0$ and the one-step transition matrix $\bm P$ given by \eqref{block-P}, we can cite Lemma~\ref{methodology} with $m=4n^2+8$ to conclude that 
\begin{align*}
f_n(p_A,p_B)&:=\P(\text{\textnormal{Team $A$ wins game to} $n$}\mid \text{\textnormal{Team $A$ serves first}}) \\
&\qquad{}-\P(\text{\textnormal{Team $A$ wins game to} $n$}\mid \text{\textnormal{Team $B$ serves first}}) \\
&\;=(\bm M \bm S)_{(0,0,2),\text{Team }A\text{ wins }}-(\bm M \bm S)_{(0,0,4),\text{Team }A\text{ wins }},
\end{align*}
where $\bm M:=(\bm I-\bm Q)^{-1}$ is the fundamental matrix.  This shows that $f_n(p_A,p_B)$ is a rational function of $p_A$ and $p_B$.

To verify symmetry,
\begin{align*}
f_n(p_A,p_B)&=\P(\text{\textnormal{Team $A$ wins game to} $n$}\mid \text{\textnormal{Team $A$ serves first}}) \\
&\qquad{}-\P(\text{\textnormal{Team $A$ wins game to} $n$}\mid \text{\textnormal{Team $B$ serves first}}) \\
&\;{}=1-\P(\text{\textnormal{Team $B$ wins game to} $n$}\mid \text{\textnormal{Team $A$ serves first}}) \\
&\qquad{}-[1-\P(\text{\textnormal{Team $B$ wins game to} $n$}\mid \text{\textnormal{Team $B$ serves first}})] \\
&\;{}=\P(\text{\textnormal{Team $B$ wins game to} $n$}\mid \text{\textnormal{Team $B$ serves first}}) \\
&\qquad{}-\P(\text{\textnormal{Team $B$ wins game to} $n$}\mid \text{\textnormal{Team $A$ serves first}}) \\
&\;{}=f_n(p_B,p_A).
\end{align*}

Finally, if $p_A=p_B=:x\in(0,1]$, then
\begin{align*}
f_n(x,x)&=\P(\text{\textnormal{Team $A$ wins game to} $n$}\mid \text{\textnormal{Team $A$ serves first}}) \\
&\qquad{}-\P(\text{\textnormal{Team $A$ wins game to} $n$}\mid \text{\textnormal{Team $B$ serves first}}) \\
&\;{}=\P(\text{\textnormal{Team $A$ wins game to} $n$}\mid \text{\textnormal{Team $A$ serves first}}) \\
&\qquad{}-[1-\P(\text{\textnormal{Team $B$ wins game to} $n$}\mid \text{\textnormal{Team $B$ serves first}})] \\
&\;{}=2\,\P(\text{\textnormal{first-serving team wins game to} $n$})-1.
\end{align*}

The reason for ruling out the case $p_A=p_B=0$ is that, if the serving team never wins a rally, the score never changes.  In this case, the two absorbing states are inaccessible and Lemma~\ref{methodology} does not apply.  
\end{proof}

\subsection{The case $n=11$}

\textit{Mathematica} yields
\begin{align*}\label{f11}
f_{11}(x,y)&=\frac{x^{11}P_{76,1430}(x,y)}{(x+y-xy)^{19}(2-x-y+xy)^{19}P_{5,11}(x,y)P_{6,15}(x,y)} \\
&\quad{}-\frac{x^{11}P_{76,1429}(x,y)}{(x+y-xy)^{19}(2-x-y+xy)^{19}P_{5,11}(x,y)P_{6,15}(x,y)} \\
&{}=-\frac{x^{11}y^{11}P_{20,121}(x,y)}{(2-x-y+xy)^{19}P_{6,15}(x,y)},
\end{align*}
where, here and below, $P_{d,t}(x,y)$ is a polynomial in $x$ and $y$ of degree $d$ with $t$ terms.  The polynomials $P_{20,121}(x,y)$, $P_{5,11}(x,y)$, and $P_{6,15}(x,y)$ are symmetric in $x$ and $y$, with the last two positive on $(0,1)\times(0,1)$.  The complete formula for $f_{11}(x,y)$ is deferred to the Appendix.
\begin{figure}[H]
\centering
\includegraphics[width=4.75in]{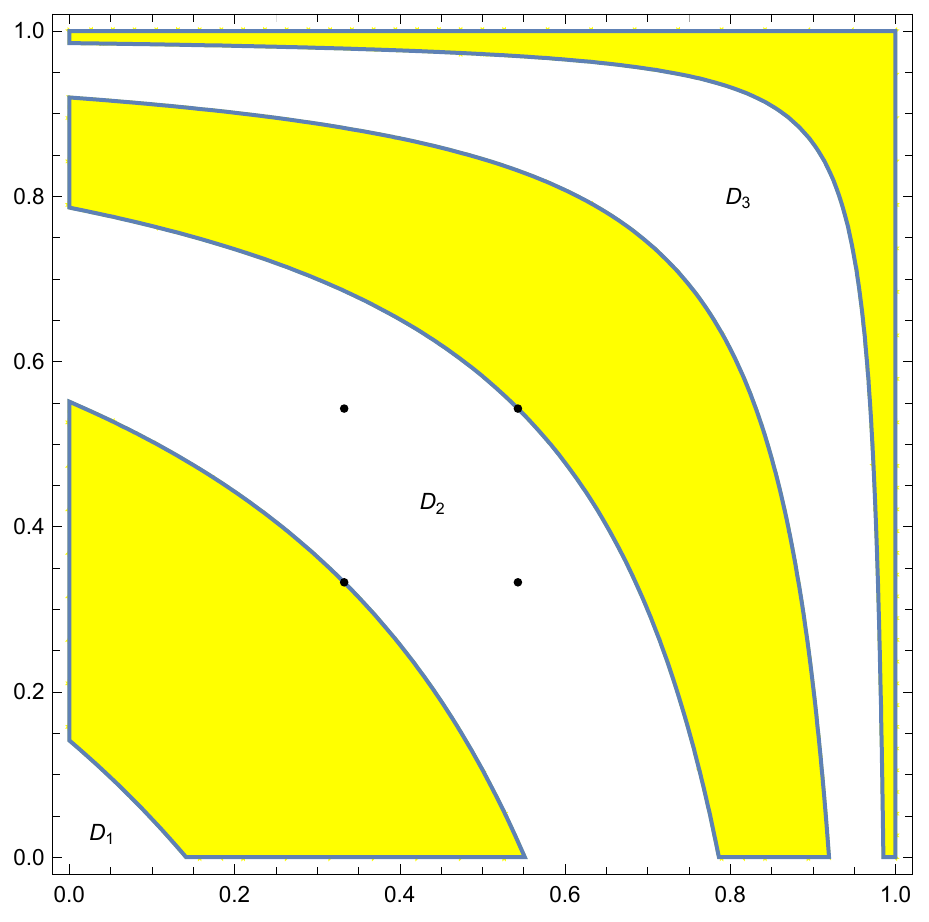}
\caption{\label{advantageplot11}The yellow region is the set of all $(p_A,p_B)$ for which the first-serving team has an advantage.  The white region is the set of all $(p_A,p_B)$ for which the first-serving team has a disadvantage.  The latter region is the union of three disjoint open connected regions, $D_1$, $D_2$, and $D_3$, each symmetric with respect to the diagonal $y=x$.  The four black dots are the corners of the square $(x_2,x_3)\times(x_2,x_3)$, where $x_2\approx0.332744$ and $x_3\approx0.543030$.  We believe that the region $D_2$ contains all $(p_A,p_B)$ that would typically be encountered in high-level pickleball.  Conventional side-out scoring to 11 is assumed.}
\end{figure}

The region in $[0,1]\times[0,1]-\{(0,0)\}$ where $f_{11}$ is negative, that is, the first-serving team has a disadvantage, is the union of three disjoint open connected regions, $D_1$, $D_2$, and $D_3$, each symmetric with respect to the diagonal $y=x$.  See Figure~\ref{advantageplot11}.  The function $f_{11}$ restricted to the diagonal $y=x$ has five zeros in $(0,1)$, namely (when rounded to six decimal places)
\begin{equation*}
0.073510,\; 0.332744,\; 0.543030,\; 0.723066,\; 0.883588.
\end{equation*}
Denoting them by $x_1$ through $x_5$, $f_{11}(x,x)$ is negative for $x$ in the intervals $(0,x_1)$, $(x_2,x_3)$, and $(x_4,x_5)$ and positive for $x$ in the intervals $(x_1,x_2)$, $(x_3,x_4)$, and $(x_5,1)$.  See Figure~\ref{winprobgraph}. 

In particular, $f_{11}(x,x)<0$ for all $x\in(x_2,x_3)$.  By virtue of Figure~\ref{advantageplot11}, $f_{11}<0$ on the square $(x_2,x_3)\times(x_2,x_3)$, which is only a small portion of $D_2$.  We believe that $D_2$ includes all pairs $(p_A,p_B)$ that would ordinarily be encountered in high-level pickleball.

We can calculate that $\min_{(x,y)\in D_2}f_{11}(x,y)\approx-7.95109\times10^{-9}$ with the minimum achieved at $(x^\circ,x^\circ)$ with $x^\circ\approx0.523681$.  This is done by minimizing $f_{11}$ over $0\le x\le0.8$, $0\le y\le0.8$, and $0.5\le x+y\le1.2$.  This closed region contains $D_2$ but does not intersect $D_1$ or $D_3$.

Thus, the answer to the question in the title for $n=11$ is, ``In a game to 11, the first-serving team has a very slight disadvantage.''  In order for this disadvantage to be detectable by simulation, the sample size $N$ (i.e., the number of simulated games) would have to be large enough that twice the standard error (about $1/\sqrt{N}$) is no larger than the largest disadvantage, about $8\times10^{-9}$.  Roughly, this would require $N$ to be of the order of $10^{16}$, which is likely beyond current computer capabilities.

\subsection{The case $n=15$}

\textit{Mathematica} yields
\begin{align*}
f_{15}(x,y)&=\frac{x^{15}P_{104,2640}(x,y)}{(x+y-xy)^{27}(2-x-y+xy)^{27}P_{5,11}(x,y)P_{6,15}(x,y)} \\
&\quad{}-\frac{x^{15}P_{104,2639}(x,y)}{(x+y-xy)^{27}(2-x-y+xy)^{27}P_{5,11}(x,y)P_{6,15}(x,y)} \\
&{}=-\frac{x^{15}y^{15}P_{28,225}(x,y)}{(2-x-y+xy)^{27}P_{6,15}(x,y)},
\end{align*}
where the polynomials $P_{28,225}(x,y)$, $P_{5,11}(x,y)$, and $P_{6,15}(x,y)$ are symmetric in $x$ and $y$, with the last two positive on $(0,1)\times(0,1)$.
The complete formula for $f_{15}(x,y)$ is deferred to the Appendix.

\begin{figure}[H]
\centering
\includegraphics[width=4.75in]{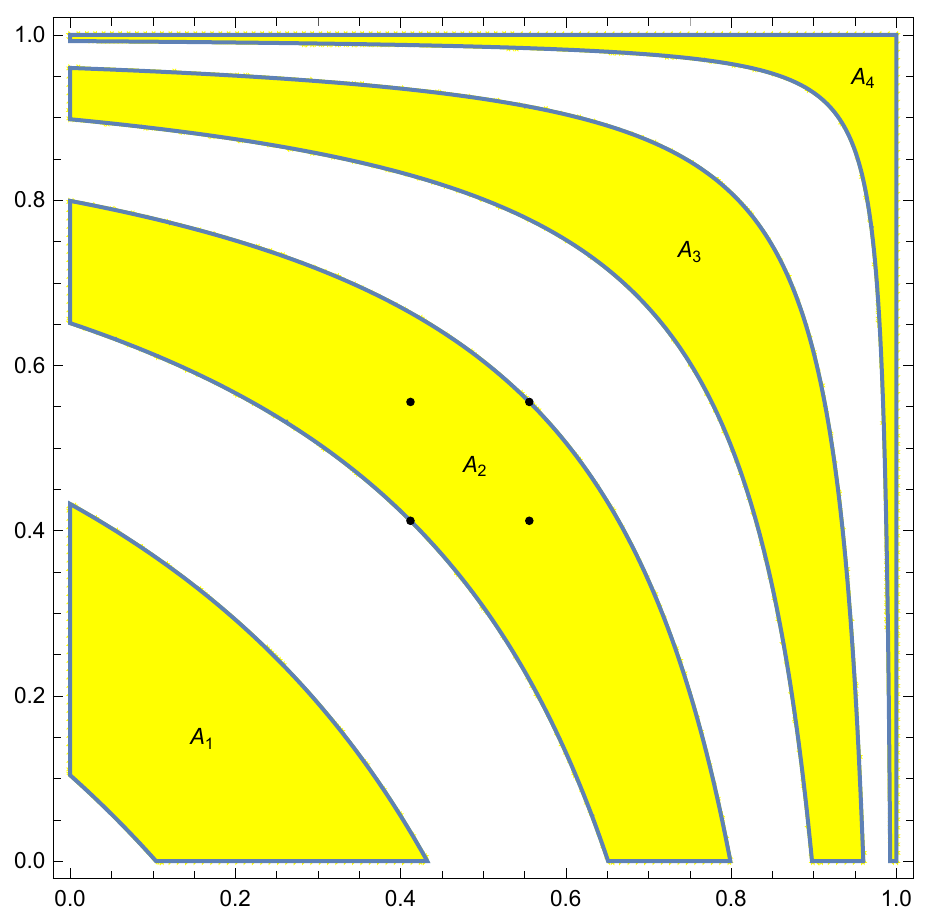}
\caption{\label{advantageplot15}The yellow region is the set of all $(p_A,p_B)$ for which the first-serving team has an advantage.  The white region is the set of all $(p_A,p_B)$ for which the first-serving team has a disadvantage.  The former region is the union of four disjoint open connected regions, $A_1$, $A_2$, $A_3$, and $A_4$, each symmetric with respect to the diagonal $y=x$.  The four black dots are the corners of the square $(x_3,x_4)\times(x_3,x_4)$, where $x_3\approx0.411942$ and $x_4\approx0.555784$.  We believe that the region $A_2$ contains all $(p_A,p_B)$ that would typically be encountered in high-level pickleball.  Conventional side-out scoring to 15 is assumed.}
\end{figure}

\begin{figure}[htb]
\centering
\includegraphics[width=4.75in]{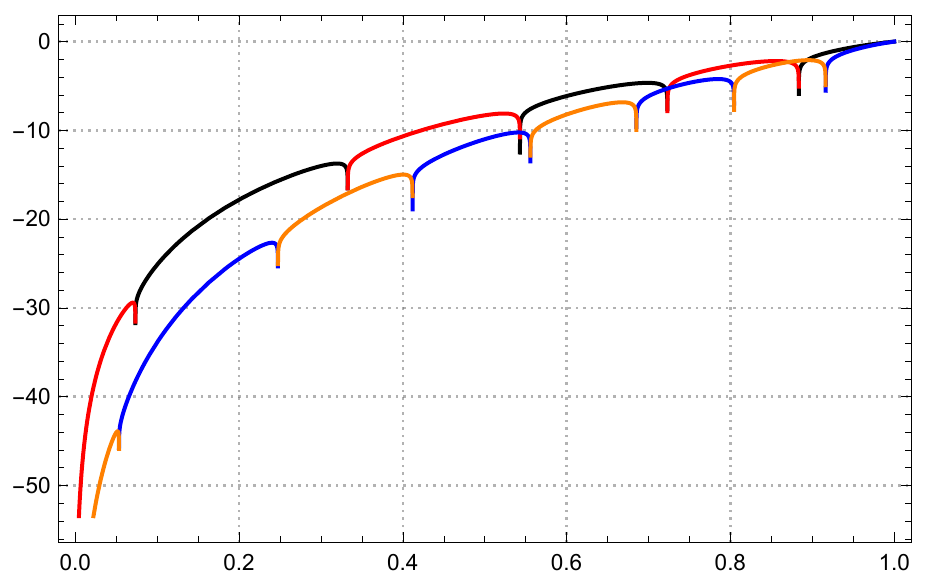}
\caption[.]{\label{winprobgraph}Graphs of $\log_{10}|f_n(x,x)|$ as a function of $x$ for $n=11$ (black and red) and for $n=15$ (blue and orange).  The black and blue (resp., red and orange) sections of the graphs correspond to intervals where the first-serving team's advantage is positive (resp., negative).  Conventional side-out scoring is assumed.}
\end{figure}

The region in $[0,1]\times[0,1]-\{(0,0)\}$ where $f_{15}$ is positive, that is, the first-serving team has an advantage, is the union of four disjoint open connected regions, $A_1$, $A_2$, $A_3$, and $A_4$, each symmetric with respect to the diagonal $y=x$.  See Figure~\ref{advantageplot15}.  The function $f_{15}$ restricted to the diagonal $y=x$ has seven zeros in $(0,1)$, namely (when rounded to six decimal places)
\begin{equation*}
0.053458,\;0.247565,\;0.411942,\;0.555784,\;0.685200,\;0.804458,\;0.916374.
\end{equation*}
Denoting them by $x_1$ through $x_7$, the function $f_{15}(x,x)$ is positive for $x$ in the intervals $(x_1,x_2)$, $(x_3,x_4)$, $(x_5,x_6)$, and $(x_7,1)$ and negative for $x$ in the intervals $(0,x_1)$, $(x_2,x_3)$, $(x_4,x_5)$, and $(x_6,x_7)$.  Again, see Figure~\ref{winprobgraph}. 

In particular, $f_{15}(x,x)>0$ for all $x\in(x_3,x_4)$.  By virtue of Figure~\ref{advantageplot15}, $f_{15}>0$ on the square $(x_3,x_4)\times(x_3,x_4)$, which is only a small portion of $A_2$.  We believe that $A_2$ includes all pairs $(p_A,p_B)$ that would ordinarily be encountered in high-level pickleball, although it may not be quite as clear as the corresponding claim for $n=11$.

We can calculate that $\max_{(x,y)\in A_2}f_{15}(x,y)\approx 5.81408\times10^{-11}$ with the maximum achieved at $(x^\circ,x^\circ)$ with $x^\circ\approx 0.541281$.  This is done by maximizing $f_{15}$ over $0\le x\le0.85$, $0\le y\le0.85$, and $0.6\le x+y\le1.2$.  This region contains $A_2$ but does not intersect $A_1$, $A_3$, or $A_4$.

Thus, the answer to the question in the title for $n=15$ is, ``In a game to 15, the first-serving team has a very slight advantage.''  This advantage is even smaller than the disadvantage when $n=11$, hence is not detectable by simulation.

\section{Modified rally scoring}

With rally scoring, a point is awarded to the winner of each rally, and that team then serves to begin the next rally.  A system called \textit{modified rally scoring} has been adopted by Major League Pickleball and is at least under consideration for adoption by USA Pickleball.  In a game to $n$, modified rally scoring coincides with rally scoring until one of the teams has reached $n-1$ points.  From there, that team can score additional points only when serving.  Furthermore, given that one team has reached $n-1$ points, the opposing team, once it reaches $n-3$ points, can score additional points only when serving.  Unlike with side-out scoring, there is no attempt to mitigate the advantage or disadvantage of the first-serving team.

A game to $n$ can be modeled by a Markov chain in a state space of $2n^2+6$ states.  $2n^2+4$ of the states are of the form $(i,j,k)$, where $(i,j)\in(\{0,1,2,\ldots,n-1\}\times\{0,1,2,\ldots,n-1\})\cup\{(n-1,n),(n,n-1)\}$ and $k\in\{1,2\}$.  In state $(i,j,k)$, Team $A$ has $i$ points, Team $B$ has $j$ points, and $k$ indicates whether it is Team $A$'s serve (1) or Team $B$'s serve (2).  There are also two absorbing states, ``Team $A$ wins'' and ``Team $B$ wins''.  The three scores $(n-1,n-1)$, $(n-1,n)$, and $(n,n-1)$, can be thought of as ``deuce'', ``advantage Team $B$'', and ``advantage Team $A$''.   The transition matrix $\bm P$ can be written in block form as in \eqref{block-P}, where $m=2n^2+4$.

Let $p_A$ (resp., $p_B$) be the probability that Team $A$ (resp., Team $B$) wins a rally when serving, let $q_A:=1-p_A$ (resp., $q_B:=1-p_B$), and define the entries of $\bm Q$ and $\bm S$ as follows.

For $i,j\in\{0,1,\ldots,n-2\}$,
\begin{small}
\begin{align*}
Q((i,j,1),(i+1,j,1))&=p_A, & Q((i,j,1),(i,j+1,2))&=q_A, \\
Q((i,j,2),(i,j+1,2))&=p_B, & Q((i,j,2),(i+1,j,1))&=q_B. 
\end{align*}
\end{small}
For $i\in\{0,1,\ldots,n-4\}$,
\begin{small}
\begin{align*}
Q((i,n-1,1),(i+1,n-1,1))&=p_A, & Q((i,n-1,1),(i,n-1,2))&=q_A, \\
S((i,n-1,2),\text{Team }B\text{ wins})&=p_B, & Q((i,n-1,2),(i+1,n-1,1))&=q_B,
\end{align*}
\end{small}
together with
\begin{small}
\begin{align*}
Q((n-3,n-1,1),(n-2,n-1,1))&=p_A, & Q((n-3,n-1,1),(n-3,n-1,2))&=q_A, \\
S((n-3,n-1,2),\text{Team }B\text{ wins})&=p_B, & Q((n-3,n-1,2),(n-3,n-1,1))&=q_B, \\
Q((n-2,n-1,1),(n-1,n-1,1))&=p_A, & Q((n-2,n-1,1),(n-2,n-1,2))&=q_A, \\
S((n-2,n-1,2),\text{Team }B\text{ wins})&=p_B, & Q((n-2,n-1,2),(n-2,n-1,1))&=q_B.
\end{align*}
\end{small}
For $j\in\{0,1,\ldots,n-4\}$,
\begin{small}
\begin{align*}
S((n-1,j,1),\text{Team }A\text{ wins})&=p_A, & Q((n-1,j,1),(n-1,j+1,2))&=q_A, \\
Q((n-1,j,2),(n-1,j+1,2))&=p_B, & Q((n-1,j,2),(n-1,j,1))&=q_B,
\end{align*}
\end{small}
together with
\begin{small}
\begin{align*}
S((n-1,n-3,1),\text{Team }A\text{ wins})&=p_A, & Q((n-1,n-3,1),(n-1,n-3,2))&=q_A, \\
Q((n-1,n-3,2),(n-1,n-2,2))&=p_B, & Q((n-1,n-3,2),(n-1,n-3,1))&=q_B, \\
S((n-1,n-2,1),\text{Team }A\text{ wins})&=p_A, & Q((n-1,n-2,1),(n-1,n-2,2))&=q_A, \\
Q((n-1,n-2,2),(n-1,n-1,2))&=p_B, & Q((n-1,n-2,2),(n-1,n-2,1))&=q_B.
\end{align*}
\end{small}
Finally,
\begin{small}
\begin{align*}
Q((n-1,n-1,1),(n,n-1,1))&=p_A, & Q((n-1,n-1,1),(n-1,n-1,2))&=q_A, \\
Q((n-1,n-1,2),(n-1,n,2))&=p_B, & Q((n-1,n-1,2),(n-1,n-1,1))&=q_B, \\
\noalign{\smallskip}
Q((n-1,n,1),(n-1,n-1,1))&=p_A, & Q((n-1,n,1),(n-1,n,2))&=q_A, \\
S((n-1,n,2),\text{Team }B\text{ wins})&=p_B, & Q((n-1,n,2),(n-1,n,1))&=q_B, \\
\noalign{\smallskip}
S((n,n-1,1),\text{Team }A\text{ wins})&=p_A, & Q((n,n-1,1),(n,n-1,2))&=q_A, \\
Q((n,n-1,2),(n-1,n-1,2))&=p_B, & Q((n,n-1,2),(n,n-1,1))&=q_B.
\end{align*}
\end{small}

\noindent All entries not specified are 0.  For specificity, it is convenient to order the states $(i,j,k)$ lexicographically, followed by states ``Team $A$ wins'' and ``Team $B$ wins''.  That completes the specification of the transition matrix $\bm P$.  The initial state is $(0,0,1)$ if Team $A$ is first server, $(0,0,2)$ if Team $B$ is first server.

\begin{theorem}\label{Thm:windiff-mr}
Consider a game of pickleball doubles with modified rally scoring.  Let $p_A$ (resp., $p_B$) be the probability that Team $A$ (resp., Team $B$) wins a rally when serving, and assume that rallies are independent and $p_A+p_B>0$. Then the conclusions of Theorem~\ref{Thm:windiff-so}, but with $f_n$ relabeled as $f_n^*$, hold in this context.
\end{theorem}

\begin{remark}
The function $f_n^*$ might be called the \textit{first-serving team's advantage} under modified rally scoring to $n$.  A negative advantage is of course a disadvantage.  \textit{Mathematica} yields
\begin{equation*}
f_{21}^*(x,y)=-\frac{xy(1-x-y)P_{56,1005}(x,y)}{(x+y-x y) P_{4,7}(x,y)},
\end{equation*}
where again $P_{d,t}(x,y)$ is a polynomial in $x$ and $y$ of degree $d$ with $t$ terms, and here both polynomials are symmetrix in $x$ and $y$ and positive on $(0,1)\times(0,1)$.  The formula for $f_{21}^*(x,y)$ would require multiple pages, so we do not include it.  But we do include a formula in the Appendix for $f_{21}^*$ restricted to the diagonal $y=x$.

It follows that $f_{21}^*(p_A,p_B)>0$, that is, the first-serving team has an advantage, if $p_A>q_B$ (i.e., if Team $A$'s probability of winning a rally is greater when Team $A$ serves than when Team $B$ serves), and that $f_{21}^*(p_A,p_B)<0$,  that is, the first-serving team has a disadvantage, if $p_A<q_B$ (i.e., if Team $A$'s probability of winning a rally is less when Team $A$ serves than when Team $B$ serves).  These conclusions seem intuitively plausible and might have been expected.

While in the case of side-out scoring the advantages for the first server are very small across the range of values of $p_A$ and $p_B$ that are likely to be encountered in high-level pickleball, here the advantages can be significant.  For example, $f_{21}^*(0.44,0.44)\approx-0.0137951$, which of course is large enough to be detectable by simulation.
\end{remark}

\begin{figure}[htb]
\centering
\includegraphics[width=4.75in]{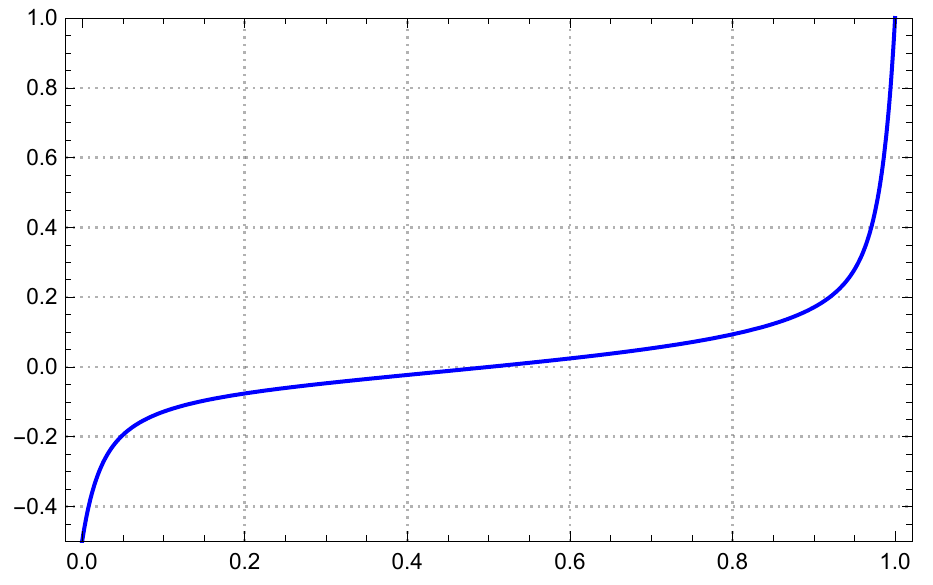}
\caption[.]{\label{windiffgraph-mr}Graph of $f_{21}^*(x,x)$ as a function of $x$.  Modified rally scoring to 21 is assumed.}
\end{figure}

\begin{proof}
The proof is virtually identical to that of Theorem~\ref{Thm:windiff-so}.

Here, the reason for ruling out the case $p_A=p_B=0$ is that, if the serving team never wins a rally, then the sequence of scores is 0-0, 0-1, 1-1, 1-2, 2-2, 2-3, 3-3, and so on, where the first-listed score is that of the first-serving team.  No team ever leads by more than one, so the game continues forever.   Again, the two absorbing states are inaccessible and Lemma~\ref{methodology} does not apply. 
\end{proof}

\section{Hybrid rally scoring: Modified rally scoring with traditional doubles server rotation}

We have not discussed server rotation, so let us consider it here.  

In side-out scoring, the first serve after a side out is always from the even court.  If the serving team wins a rally, the two players on that team switch places and the player who served continues to serve.  In no other cases do the players on a team switch places.  If a team's score is even (resp., odd), then the player on that team who began the game on the even court is on the even (resp., odd) court.

In modified rally scoring, there is no switching places.  A player who begins the game on the even court remains there throughout the game.  When the serving team's score is even (resp., odd), it serves from the even (resp., odd) court.

We would like to propose a hybrid system, comprising modified rally scoring and traditional doubles server rotation.  This has the advantages of rally scoring while being more familiar to players used to side-out scoring.

As in side-out scoring, each team serves until it has faulted twice, except at the beginning of the game, in which case the first-serving team serves until it has faulted once.  Also, the first serve after a side out is always from the even court.  If the serving team wins a rally, the two players on that team switch places and the player who served continues to serve.  If the receiving team wins a rally, the two players on that team switch places.  (This differs from traditional doubles.  It is needed to ensure that player positions align with a team's score.)  If a team's score is even (resp., odd), then the player on that team who began the game on the even court is on the even (resp., odd) court.

The state space for the Markov chain that models a game under hybrid rally scoring is exactly as for side-out scoring.  Let $p_A$ (resp., $p_B$) be the probability that Team $A$ (resp., Team $B$) wins a rally when serving, let $q_A:=1-p_A$ (resp., $q_B:=1-p_B$), and define the entries of $\bm Q$ and $\bm S$ as follows.

For $i,j\in\{0,1,\ldots,n-2\}$,
\begin{small}
\begin{align*}
Q((i,j,1),(i+1,j,1))&=p_A, & Q((i,j,1),(i,j+1,2))&=q_A, \\
Q((i,j,2),(i+1,j,2))&=p_A, & Q((i,j,2),(i,j+1,3))&=q_A, \\
Q((i,j,3),(i,j+1,3))&=p_B, & Q((i,j,3),(i+1,j,4))&=q_B, \\
Q((i,j,4),(i,j+1,4))&=p_B, & Q((i,j,4),(i+1,j,1))&=q_B.
\end{align*}
\end{small}
For $i\in\{0,1,\ldots,n-4\}$,
\begin{small}
\begin{align*}
Q((i,n-1,1),(i+1,n-1,1))&=p_A, & Q((i,n-1,1),(i,n-1,2))&=q_A, \\
Q((i,n-1,2),(i+1,n-1,2))&=p_A, & Q((i,n-1,2),(i,n-1,3))&=q_A, \\
S((i,n-1,3),\text{Team }B\text{ wins})&=p_B, & Q((i,n-1,3),(i+1,n-1,4))&=q_B, \\
S((i,n-1,4),\text{Team }B\text{ wins})&=p_B, & Q((i,n-1,4),(i+1,n-1,1))&=q_B, 
\end{align*}
\end{small}
together with
\begin{small}
\begin{align*}
Q((n-3,n-1,1),(n-2,n-1,1))&=p_A, & Q((n-3,n-1,1),(n-3,n-1,2))&=q_A, \\
Q((n-3,n-1,2),(n-2,n-1,2))&=p_A, & Q((n-3,n-1,2),(n-3,n-1,3))&=q_A, \\
S((n-3,n-1,3),\text{Team }B\text{ wins})&=p_B, & Q((n-3,n-1,3),(n-3,n-1,4))&=q_B, \\
S((n-3,n-1,4),\text{Team }B\text{ wins})&=p_B, & Q((n-3,n-1,4),(n-3,n-1,1))&=q_B,\\
\noalign{\smallskip}
Q((n-2,n-1,1),(n-1,n-1,1))&=p_A, & Q((n-2,n-1,1),(n-2,n-1,2))&=q_A, \\
Q((n-2,n-1,2),(n-1,n-1,2))&=p_A, & Q((n-2,n-1,2),(n-2,n-1,3))&=q_A, \\
S((n-2,n-1,3),\text{Team }B\text{ wins})&=p_B, & Q((n-2,n-1,3),(n-2,n-1,4))&=q_B, \\
S((n-2,n-1,4),\text{Team }B\text{ wins})&=p_B, & Q((n-2,n-1,4),(n-2,n-1,1))&=q_B.
\end{align*}
\end{small}
For $j\in\{0,1,\ldots,n-4\}$,
\begin{small}
\begin{align*}
S((n-1,j,1),\text{Team }A\text{ wins})&=p_A, & Q((n-1,j,1),(n-1,j+1,2))&=q_A, \\
S((n-1,j,2),\text{Team }A\text{ wins})&=p_A, & Q((n-1,j,2),(n-1,j+1,3))&=q_A, \\
Q((n-1,j,3),(n-1,j+1,3))&=p_B, & Q((n-1,j,3),(n-1,j,4))&=q_B, \\
Q((n-1,j,4),(n-1,j+1,4))&=p_B, & Q((n-1,j,4),(n-1,j,1))&=q_B,
\end{align*}
\end{small}
together with
\begin{small}
\begin{align*}
S((n-1,n-3,1),\text{Team }A\text{ wins})&=p_A, & Q((n-1,n-3,1),(n-1,n-3,2))&=q_A, \\
S((n-1,n-3,2),\text{Team }A\text{ wins})&=p_A, & Q((n-1,n-3,2),(n-1,n-3,3))&=q_A, \\
Q((n-1,n-3,3),(n-1,n-2,3))&=p_B, & Q((n-1,n-3,3),(n-1,n-3,4))&=q_B, \\
Q((n-1,n-3,4),(n-1,n-2,4))&=p_B, & Q((n-1,n-3,4),(n-1,n-3,1))&=q_B,\\
\noalign{\smallskip}
S((n-1,n-2,1),\text{Team }A\text{ wins})&=p_A, & Q((n-1,n-2,1),(n-1,n-2,2))&=q_A, \\
S((n-1,n-2,2),\text{Team }A\text{ wins})&=p_A, & Q((n-1,n-2,2),(n-1,n-2,3))&=q_A, \\
Q((n-1,n-2,3),(n-1,n-1,3))&=p_B, & Q((n-1,n-2,3),(n-1,n-2,4))&=q_B, \\
Q((n-1,n-2,4),(n-1,n-1,4))&=p_B, & Q((n-1,n-2,4),(n-1,n-2,1))&=q_B.
\end{align*}
\end{small}
Finally,
\begin{small}
\begin{align*}
Q((n-1,n-1,1),(n,n-1,1))&=p_A, & Q((n-1,n-1,1),(n-1,n-1,2))&=q_A, \\
Q((n-1,n-1,2),(n,n-1,2))&=p_A, & Q((n-1,n-1,2),(n-1,n-1,3))&=q_A, \\
Q((n-1,n-1,3),(n-1,n,3))&=p_B, & Q((n-1,n-1,3),(n-1,n-1,4))&=q_B, \\
Q((n-1,n-1,4),(n-1,n,4))&=p_B, & Q((n-1,n-1,4),(n-1,n-1,1))&=q_B, \\
\noalign{\medskip}
Q((n-1,n,1),(n-1,n-1,1))&=p_A, & Q((n-1,n,1),(n-1,n,2))&=q_A, \\
Q((n-1,n,2),(n-1,n-1,2))&=p_A, & Q((n-1,n,2),(n-1,n,3))&=q_A, \\
S((n-1,n,3),\text{Team }B\text{ wins})&=p_B, & Q((n-1,n,3),(n-1,n,4))&=q_B, \\
S((n-1,n,4),\text{Team }B\text{ wins})&=p_B, & Q((n-1,n,4),(n-1,n,1))&=q_B, \\ 
\noalign{\medskip}
S((n,n-1,1),\text{Team }A\text{ wins})&=p_A, & Q((n,n-1,1),(n,n-1,2))&=q_A, \\
S((n,n-1,2),\text{Team }A\text{ wins})&=p_A, & Q((n,n-1,2),(n,n-1,3))&=q_A, \\
Q((n,n-1,3),(n-1,n-1,3))&=p_B, & Q((n,n-1,3),(n,n-1,4))&=q_B, \\
Q((n,n-1,4),(n-1,n-1,4))&=p_B, & Q((n,n-1,4),(n,n-1,1))&=q_B.
\end{align*}
\end{small}

All entries not specified are 0.  For specificity, it is convenient to order the states $(i,j,k)$ lexicographically, followed by states ``Team $A$ wins'' and ``Team $B$ wins''.  That completes the specification of the transition matrix $\bm P$.  The initial state is $(0,0,2)$ if Team $A$ is first server, $(0,0,4)$ if Team $B$ is first server.

\begin{theorem}\label{Thm:windiff-hr}
Consider a game of pickleball doubles with hybrid rally scoring.  Let $p_A$ (resp., $p_B$) be the probability that Team $A$ (resp., Team $B$) wins a rally when serving. Assume that rallies are independent and $p_A+p_B>0$. Then the conclusions of Theorem~\ref{Thm:windiff-so}, but with $f_n$ relabeled as $f_n^\circ$, hold in this context.
\end{theorem}

\begin{figure}[htb]
\centering
\includegraphics[width=4.75in]{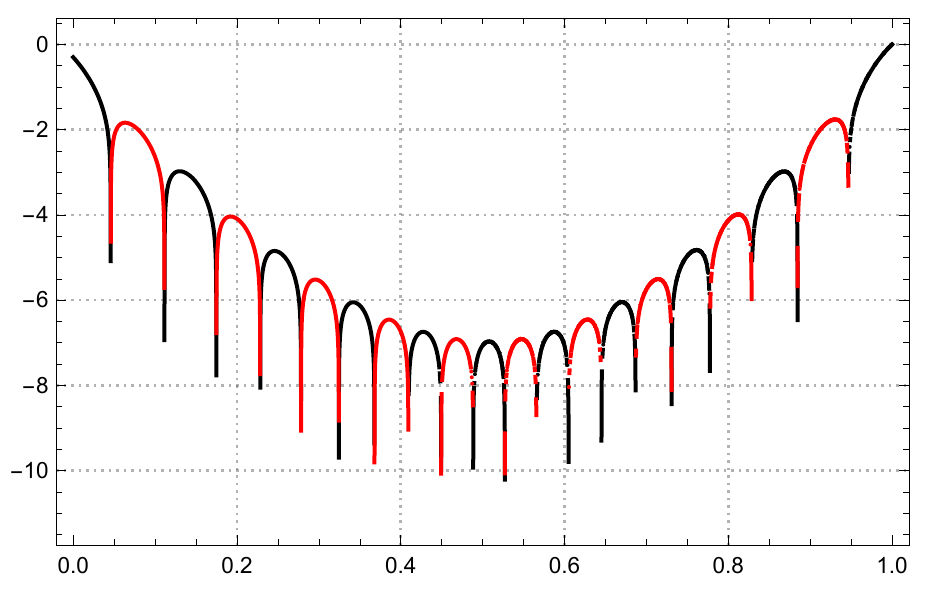}
\caption[.]{\label{winprobgraph-hr}Graph of $\log_{10}|f_{21}^\circ(x,x)|$ as a function of $x$.  The black (resp., red) sections of the graphs correspond to intervals where the first-serving team's advantage is positive (resp., negative).  Hybrid rally scoring to 21 is assumed.}
\end{figure}

\begin{remark}
The function $f_n^\circ$ might be called the \textit{first-serving team's advantage} under hybrid rally scoring to $n$.  A negative advantage is of course a disadvantage.  \textit{Mathematica} yields
\begin{equation*}
f_{21}^\circ(x,y)=\frac{xy\,P_{67,1459}(x,y)}{(x+y-xy)(2-x-y+xy)P_{5,11}(x,y)P_{6,15}(x,y)},
\end{equation*}
where each polynomial is symmetric in $x$ and $y$, and the denominator polynomials are in addition positive on $(0,1)\times(0,1)$.  The formula for $f_{21}^\circ(x,y)$ would require multiple pages, so we do not include it.  But we do include a formula in the Appendix for $f_{21}^\circ$ restricted to the diagonal $y=x$.  We were unable to produce a plot analogous to Figures~\ref{advantageplot11} and \ref{advantageplot15}.

The rational function $f_{21}^\circ(x,x)$ of $x$ has 20 zeros in $(0,1)$, namely (when rounded to six decimal places)
\begin{align*}
& 0.045911,\; 0.111525,\; 0.174871,\; 0.228450,\; 0.278268, \\
& 0.324458,\; 0.367935,\; 0.409386,\; 0.449396,\; 0.488480, \\
& 0.527122,\; 0.565789,\; 0.604954,\; 0.645110,\; 0.686804, \\
& 0.730699,\; 0.777620,\; 0.828561,\; 0.884621,\; 0.946746.
\end{align*}
Denote them by $x_1$ through $x_{20}$ and put $x_0:=0$ and $x_{21}:=1$.  Then $f_{21}^*(x,x)$ is positive for $x$ in the interval $(x_{2n},x_{2n+1})$ for $n=0,1,\ldots,10$ and negative for $x$ in the interval $(x_{2n-1},x_{2n})$ for $n=1,2,\ldots,10$.  See Figure~\ref{winprobgraph-hr}.  

Notice that $x_8$, $x_9$, and $x_{10}$ belong to the interval $[0.4,0.5]$.
\end{remark}

\begin{proof}
The proof is virtually identical to that of Theorem~\ref{Thm:windiff-so}.

Here, the reason for ruling out the case $p_A=p_B=0$ is that, if the serving team never wins a rally, then the sequence of scores is 0-0, 0-1, 1-1, 2-1, 2-2, 2-3, 3-3, 4-3, and so on, where the first-listed score is that of the first-serving team.  No team ever leads by more than one, so the game continues forever.  Again, the two absorbing states are inaccessible and Lemma~\ref{methodology} does not apply. 
\end{proof}

\section{Side-out scoring vs.\ modified rally scoring}

Here we compare side-out scoring to 11 with modified rally scoring to 21.  We consider the probability that Team $A$ wins, the expected duration of the game (duration = number of rallies), and the standard deviation of the game's duration.  In each case we assume that the first-serving team is chosen by the toss of a fair coin (or an equivalent randomization).  Results could be expressed in formulas, tables, spreadsheets, or graphics.  We prefer the graphical method used by Newton and Keller (2005), namely simultaneous graphs of evenly spaces cross-sections of the function of two variables.

Let
\begin{align*}
P(p_A,p_B)&:=\P(\text{Team $A$ wins game to 11 under side-out scoring}), \\
P^*(p_A,p_B)&:=\P(\text{Team $A$ wins game to 21 under modified rally scoring}).
\end{align*}
Each probability is the average of the two corresponding probabilities, one with Team $A$ as first server, the other with Team $B$ as first server.  
\begin{figure}[H]
\centering
\includegraphics[width=4.75in]{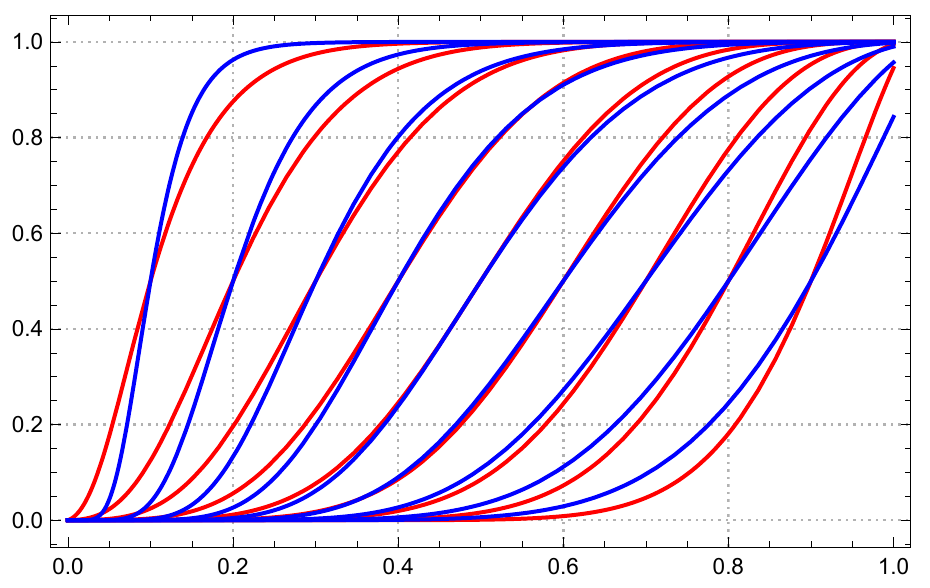}
\caption{\label{winprobgraphs}Graphs of the probability that Team $A$ wins a game under side-out scoring to 11 (blue) and under modified rally scoring to 21 (red) as a function of $p_A$, when $p_B=m/10$ ($m=1,2,\ldots,9$).  We assume that the first server is determined by the toss of a fair coin.  To distinguish graphs of the same color, note that the graph for $p_B=m/10$ has win probability approximately equal to 1/2 at $p_A=m/10$.}
\end{figure}

In Figure~\ref{winprobgraphs} we graph $P(\cdot,p_B)$ and $P^*(\cdot,p_B)$ for $p_B=m/10$ ($m=1,2,\ldots,9$).  Especially for $p_B=0.4, 0.5, 0.6$ there is not much difference between the two scoring systems as far as the probability that Team $A$ wins is concerned.  Having computed the probabilities that Team $A$ wins a game with Team $A$ as first server and with Team $B$ as first server, it is straightforward to extract the probability that Team $A$ wins a match consisting of multiple games.

\newpage

\begin{figure}[htb]
\centering
\includegraphics[width=4.75in]{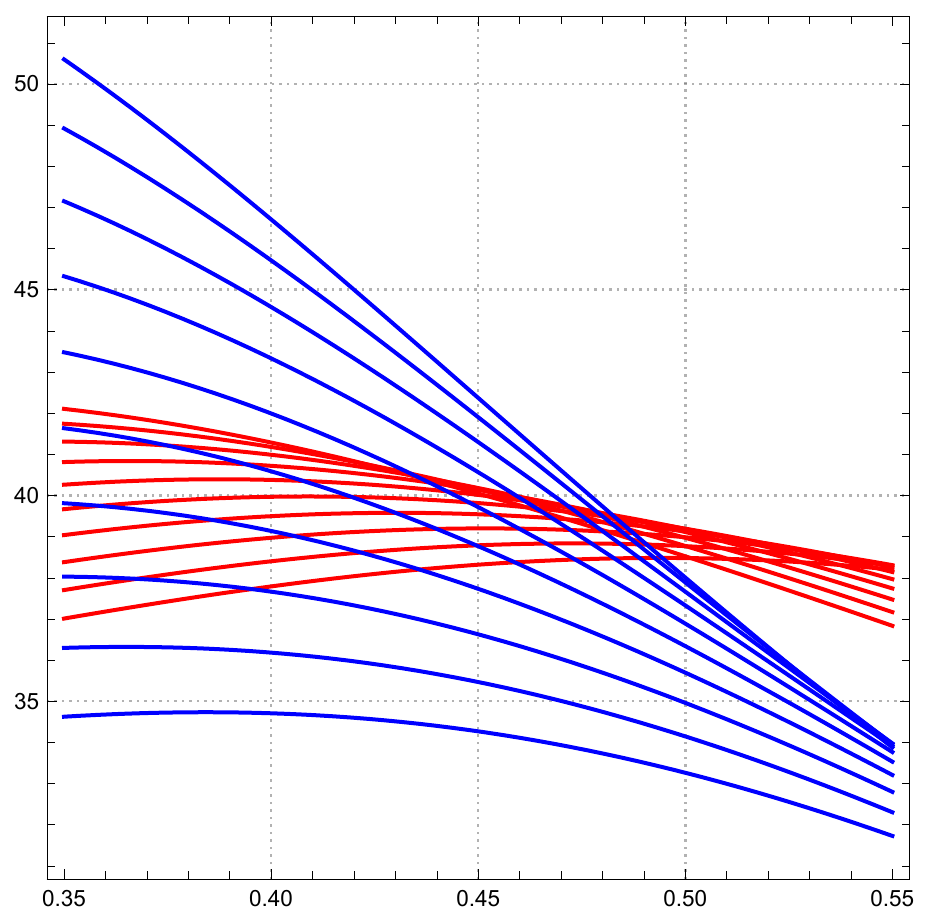}
\caption{\label{Edurgraphs}Graphs of the expected duration of a game under side-out scoring to 11 (blue) and under modified rally scoring to 21 (red) as a function of $p_A$, when $p_B=0.36,0.38,0.40,0.42,0.44,0.46,0.48,0.50,0.52,0.54$.  We assume that the first server is determined by the toss of a fair coin.  To distinguish graphs of the same color, note that, at $p_A=0.35$ the expected durations are decreasing as $p_B$ increases.}
\end{figure}

Let
\begin{align*}
E(p_A,p_B)&:=\E[\text{duration of game to 11 under side-out scoring}], \\
E^*(p_A,p_B)&:=\E[\text{duration of game to 21 under modified rally scoring}].
\end{align*}

\noindent Each expectation is the average of the two corresponding expectations, one with Team $A$ as first server, the other with Team $B$ as first server.
In Figure~\ref{Edurgraphs} we graph $E(\cdot,p_B)$ and $E^*(\cdot,p_B)$ for $p_B=0.36, 0.38, 0.40, 0.42, 0.44, 0.46, 0.48, \linebreak0.50, 0.52, 0.54$ and $0.35\le p_A\le0.55$.  We see that, at least for this range of typical values of $p_A$ and $p_B$, there is not much difference in expected game duration under side-out scoring to 11 and under modified rally scoring to 21.  Furthermore, it is not difficult to derive the expected duration of a match consisting of multiple games.

\begin{figure}[htb]
\centering
\includegraphics[width=4.75in]{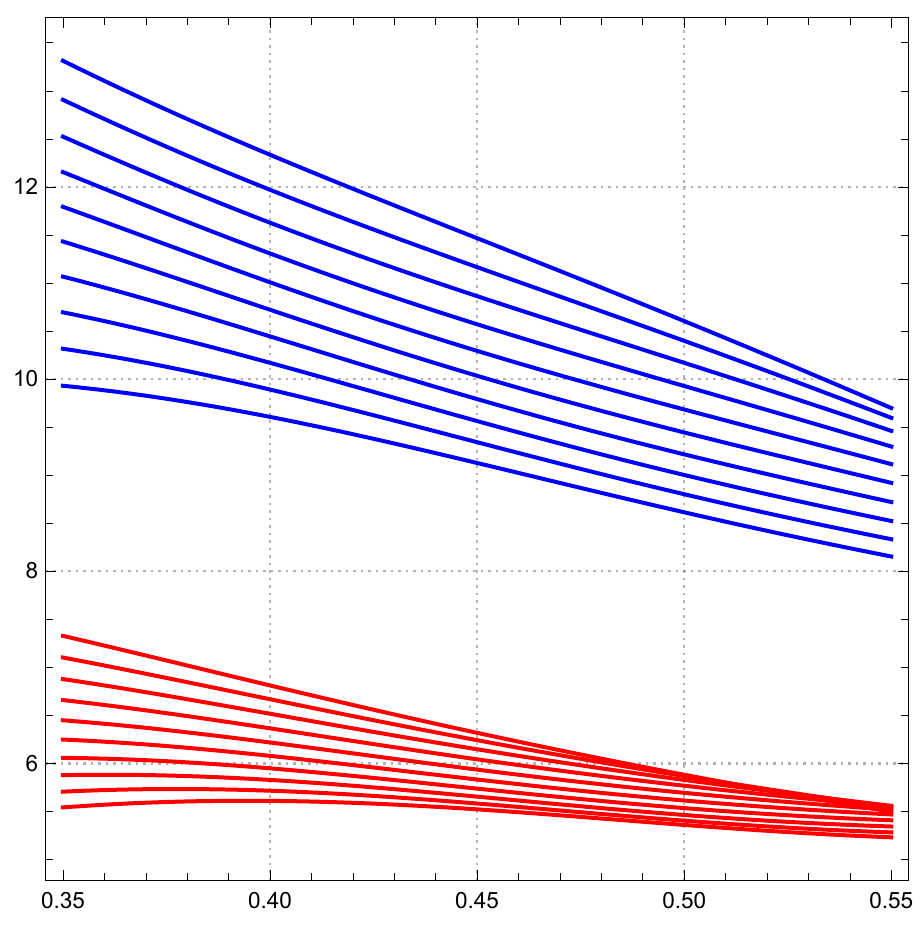}
\caption{\label{SDdurgraphs}Graphs of the standard deviation of the duration of a game under side-out scoring to 11 (blue) and under modified rally scoring to 21 (red) as a function of $p_A$, when $p_B=0.36,0.38,0.40,0.42,0.44,0.46,0.48,0.50,0.52,0.54$.  We assume that the first server is determined by the toss of a fair coin.  To distinguish graphs of the same color, note that, at $p_A=0.35$ the standard deviations of the duration are decreasing as $p_B$ increases.}
\end{figure}

Next, let
\begin{align*}
S\!D(p_A,p_B)&:=\SD(\text{duration of game to 11 under side-out scoring}), \\
S\!D^*(p_A,p_B)&:=\SD(\text{duration of game to 21 under modified rally scoring}).
\end{align*}
Each standard deviation is the square root of the variance (or of the second moment minus the square of the mean).  Each second moment is the average of the two corresponding second moments, one with Team $A$ as first server, the other with Team $B$ as first server.  Similarly, each mean is the average of the two corresponding means, one with Team $A$ as first server, the other with Team $B$ as first server.  In Figure~\ref{SDdurgraphs} we graph $S\!D(\cdot,p_B)$ and $S\!D^*(\cdot,p_B)$ for $p_B=0.36, 0.38, 0.40, 0.42, 0.44, 0.46, 0.48, 0.50, 0.52, 0.54$ and $0.35\le p_A\le0.55$.  Here the standard deviation of game duration, at least for this range of typical values of $p_A$ and $p_B$, is substantially smaller under modified rally scoring to 21 than under side-out scoring to 11.  Indeed, this is perhaps the primary benefit of modified rally scoring.  Furthermore, it is not difficult to derive the standard deviation of the duration of a match consisting of multiple games.

There are several other statistics that could be compared.  For example, the expected duration of a game could be decomposed into the expected number of serves by Team $A$ and the expected number of serves by Team $B$.  The probability of a comeback from a specified deficit could be evaluated.  We leave these issues to the interested reader.

\section{Conclusions}

Our absorbing Markov chain approach provides an answer to the title question for the most commonly used scoring methods in tournament play, showing that for side-out scoring the sign of the first-serving team's advantage depends on the number of points needed to win, while for modified rally scoring it depends on the sign of $p_A-q_B$. This method also allowed us to quantify the reduction in standard deviation of game duration that is one of the benefits of modified rally scoring.  We have also proposed and analyzed a new hybrid scoring method that maintains the lower standard deviation in game duration of modified rally scoring while also maintaining the familiar doubles server rotation from side-out scoring.

\section*{Acknowledgments}

We thank David Beecher for creating Figure~\ref{fig:pickleball-deuce}.

\begin{newreferences}

\item Calhoun, William, Dargahi-Noubary, G. R., and Shi, Yixun (2002) Volleyball scoring systems.  \textit{Mathematics and Computer Education} \textbf{36} (1) 70--79.

\item Golden, Jessica (2023) Pickleball popularity exploded last year, with more than 36 million playing the sport.  \textit{CNBC}, Jan.\ 4. \url{https://www.cnbc.com/2023/01/05/pickleball-popularity-explodes-with-more-than-36-million-playing.html}

\item Kemeny, John G. and Snell, J. Laurie (1976) \textit{Finite Markov Chains}. Springer-Verlag, New York.

\item MacPhee, I. M., Rougier, Jonathan, and Pollard, G. H. (2004) Server advantage in tennis matches.  \textit{J. Applied Probability} \textbf{41} (4) 1182--1186.

\item Newton, Paul K. and Keller, Joseph B. (2005) Probability of winning at tennis I. Theory and data.  \textit{Studies in Applied Mathematics} \textbf{114} 241--269.

\item Newton, P. K. and Pollard, G. H. (2004) Service neutral scoring strategies in tennis.  \textit{Proceedings of the Seventh Australasian Conference on Mathematics \& Computers in Sport}. (Ed.\ R. Hugh Morton and S. Ganesalingam.)  Massey University, pp.\ 221--225.


\item Xu, William (2022) The definitive mathematical solution to ``side or serve?'' Advanced pickleball strategy blog (December 19). \url{https://advancedpickleballstrategy.com/2022/12/19/the-definitive-mathematical-solution-to-side-or-serve/}

\end{newreferences}

\section*{Appendix}
The formula for $f_{11}$ in Theorem~\ref{Thm:windiff-so} is
\begin{small}
\begin{align*}
&f_{11}(x,y) \nonumber\\
&{}=-x^{11} y^{11}(252 - 2646\,x - 2646\,y + 6300\,x^2 + 15372\,x y + 6300\,y^2 + 840\,x^3 \\
&\quad{}- 10710\,x^2 y - 10710\,x y^2 + 840\,y^3 - 18060\,x^4 - 74760\,x^3 y - 106470\,x^2 y^2 \\
&\quad{}- 74760\,x y^3 - 18060\,y^4 + 17766\,x^5 + 164850\,x^4 y + 408240\,x^3 y^2 + 408240\,x^2 y^3 \\
&\quad{}+ 164850\,x y^4 + 17766\,y^5 - 588\,x^6 - 98028\,x^5 y - 506520\,x^4 y^2 - 819000\,x^3 y^3 \\
&\quad{}- 506520\,x^2 y^4 - 98028\,x y^5 - 588\,y^6 - 4980\,x^7 - 30198\,x^6 y + 119070\,x^5 y^2 \\
&\quad{}+ 556920\,x^4 y^3 + 556920\,x^3 y^4 + 119070\,x^2 y^5 - 30198\,x y^6 - 4980\,y^7 + 810\,x^8 \\
&\quad{}+ 43608\,x^7 y + 230958\,x^6 y^2 + 332640\,x^5 y^3 + 267120\,x^4 y^4 + 332640\,x^3 y^5 \\
&\quad{}+ 230958\,x^2 y^6 + 43608\,x y^7 + 810\,y^8 + 295\,x^9 - 4581\,x^8 y - 158220\,x^7 y^2 \\
&\quad{}- 695016\,x^6 y^3 - 1217160\,x^5 y^4 - 1217160\,x^4 y^5 - 695016\,x^3 y^6 - 158220\,x^2 y^7 \\
&\quad{}- 4581\,x y^8 + 295\,y^9 + 10\,x^{10} - 2798\,x^9 y + 5238\,x^8 y^2 + 296208\,x^7 y^3 + 1072512\,x^6 y^4 \\
&\quad{}+ 1545264\,x^5 y^5 + 1072512\,x^4 y^6 + 296208\,x^3 y^7 + 5238\,x^2 y^8 - 2798\,x y^9 + 10\,y^{10} \\
&\quad{}- 109\,x^{10} y + 11583\,x^9 y^2 + 22788\,x^8 y^3 - 271152\,x^7 y^4 - 836136\,x^6 y^5 - 836136\,x^5 y^6 \\
&\quad{}- 271152\,x^4 y^7 + 22788\,x^3 y^8 + 11583\,x^2 y^9 - 109\,x y^{10} + 531\,x^{10} y^2- 27336\,x^9 y^3 \\
&\quad{}- 91476\,x^8 y^4 + 27216\,x^7 y^5 + 179928\,x^6 y^6 + 27216\,x^5 y^7 - 91476\,x^4 y^8 - 27336\,x^3 y^9 \\
&\quad{}+ 531\,x^2 y^{10} - 1524\,x^{10} y^3 + 40110\,x^9 y^4 + 150066\,x^8 y^5 + 190512\,x^7 y^6 + 190512\,x^6 y^7 \\
&\quad{}+ 150066\,x^5 y^8 + 40110\,x^4 y^9 - 1524\,x^3 y^{10} + 2856\,x^{10} y^4 - 37044\,x^9 y^5 - 135324\,x^8 y^6 \\
&\quad{}- 183600\,x^7 y^7 - 135324\,x^6 y^8 - 37044\,x^5 y^9 + 2856\,x^4 y^{10} - 3654\,x^{10} y^5 + 20118\,x^9 y^6 \\
&\quad{}+ 66852\,x^8 y^7 + 66852\,x^7 y^8 + 20118\,x^6 y^9 - 3654\,x^5 y^{10} + 3234\,x^{10} y^6 - 4488\,x^9 y^7 \\
&\quad{}- 13878\,x^8 y^8 - 4488\,x^7 y^9 + 3234\,x^6 y^{10} - 1956\,x^{10} y^7 - 1269\,x^9 y^8 - 1269\,x^8 y^9 \\
&\quad{}- 1956\,x^7 y^{10} + 774\,x^{10} y^8 + 1010\,x^9 y^9 + 774\,x^8 y^{10} - 181\,x^{10} y^9 - 181\,x^9 y^{10} \\
&\quad{}+ 19\,x^{10} y^{10}) \\
&{}/[(2 - x - y + x y)^{19}(4\,x + 4\,y - 4\,x^2 - 10\,x y - 4\,y^2 + x^3 + 10\,x^2 y + 10\,x y^2 + y^3 \\
&\quad{}- 4\,x^3 y - 11\,x^2 y^2 - 4\,x y^3 + 5\,x^3 y^2 + 5\,x^2 y^3 - 2\,x^3 y^3)].
\end{align*}
\end{small}
\!\!Notice that the symmetry of $f_{11}$ can be verified by inspection by checking that the coefficient of $x^i y^j$ is equal to that of $x^j y^i$, in both numerator and denominator.

The formula for $f_{15}$ in Theorem~\ref{Thm:windiff-so} is
\begin{small}
\begin{align*}
&f_{15}(x,y) \\
&{}=-x^{15} y^{15} (3432 - 49764\,x - 49764\,y + 168168\,x^2 + 387816\,x y + 168168\,y^2 + 24024\,x^3 \\
&\quad{}- 276276\,x^2 y - 276276\,x y^2 + 24024\,y^3 - 1069068\,x^4 - 4348344\,x^3 y - 6378372\,x^2 y^2 \\
&\quad{}- 4348344\,x y^3 - 1069068\,y^4 + 1879878\,x^5 + 13831818\,x^4 y + 32168136\,x^3 y^2 \\
&\quad{}+ 32168136\,x^2 y^3 + 13831818\,x y^4 + 1879878\,y^5 - 636636\,x^6 - 13609596\,x^5 y \\
&\quad{}- 55591536\,x^4 y^2 - 85261176\,x^3 y^3 - 55591536\,x^2 y^4 - 13609596\,x y^5 - 636636\,y^6 \\
&\quad{}- 1232088\,x^7 - 4570566\,x^6 y + 14756742\,x^5 y^2 + 63243180\,x^4 y^3 + 63243180\,x^3 y^4 \\
&\quad{}+ 14756742\,x^2 y^5 - 4570566\,x y^6 - 1232088\,y^7 + 1108536\,x^8 + 17805216\,x^7 y \\
&\quad{}+ 73939866\,x^6 y^2 + 132107976\,x^5 y^3 + 148504356\,x^4 y^4 + 132107976\,x^3 y^5 \\
&\quad{}+ 73939866\,x^2 y^6 + 17805216\,x y^7 + 1108536\,y^8 - 12012\,x^9 - 9422556\,x^8 y \\
&\quad{}- 93381288\,x^7 y^2 - 320275956\,x^6 y^3 - 548810262\,x^5 y^4 - 548810262\,x^4 y^5 \\
&\quad{}- 320275956\,x^3 y^6 - 93381288\,x^2 y^7 - 9422556\,x y^8 - 12012\,y^9 - 216216\,x^{10} \\
&\quad{}- 1987128\,x^9 y + 25832664\,x^8 y^2 + 231735504\,x^7 y^3 + 656167512\,x^6 y^4 \\
&\quad{}+ 901800900\,x^5 y^5 + 656167512\,x^4 y^6 + 231735504\,x^3 y^7 + 25832664\,x^2 y^8 \\
&\quad{}- 1987128\,x y^9 - 216216\,y^{10} + 17472\,x^{11} + 2465892\,x^{10} y + 20456436\,x^9 y^2 \\
&\quad{}+ 6513936\,x^8 y^3 - 230173944\,x^7 y^4 - 596906310\,x^6 y^5 - 596906310\,x^5 y^6 \\
&\quad{}- 230173944\,x^4 y^7 + 6513936\,x^3 y^8 + 20456436\,x^2 y^9 + 2465892\,x y^{10} + 17472\,y^{11} \\
&\quad{}+ 13286\,x^{12} - 55068\,x^{11} y - 12001704\,x^{10} y^2 - 89715340\,x^9 y^3 - 209193270\,x^8 y^4 \\
&\quad{}- 188948760\,x^7 y^5 - 115945830\,x^6 y^6 - 188948760\,x^5 y^7 - 209193270\,x^4 y^8 \\
&\quad{}- 89715340\,x^3 y^9 - 12001704\,x^2 y^{10} - 55068\,x y^{11} + 13286\,y^{12} + 973\,x^{13} \\
&\quad{}- 157937\,x^{12} y - 616512\,x^{11} y^2 + 31513768\,x^{10} y^3 + 219642709\,x^9 y^4 \\
&\quad{}+ 569606895\,x^8 y^5 + 824503680\,x^7 y^6 + 824503680\,x^6 y^7 + 569606895\,x^5 y^8 \\
&\quad{}+ 219642709\,x^4 y^9 + 31513768\,x^3 y^{10} - 616512\,x^2 y^{11} - 157937\,x y^{12} + 973\,y^{13} \\
&\quad{}+ 14\,x^{14} - 13298\,x^{13} y + 838916\,x^{12} y^2 + 5223452\,x^{11} y^3 - 43842656\,x^{10} y^4 \\
&\quad{}- 317981950\,x^9 y^5 - 774812610\,x^8 y^6 - 1008544680\,x^7 y^7 - 774812610\,x^6 y^8 \\
&\quad{}- 317981950\,x^5 y^9 - 43842656\,x^4 y^{10} + 5223452\,x^3 y^{11} + 838916\,x^2 y^{12} - 13298\,x y^{13} \\
&\quad{}+ 14\,y^{14} - 209\,x^{14} y + 83317\,x^{13} y^2 - 2606318\,x^{12} y^3 - 18712694\,x^{11} y^4 \\
&\quad{}+ 14375647\,x^{10} y^5 + 246308205\,x^9 y^6 + 548442180\,x^8 y^7 + 548442180\,x^7 y^8 \\
&\quad{}+ 246308205\,x^6 y^9 +14375647\,x^5 y^{10} - 18712694\,x^4 y^{11} - 2606318\,x^3 y^{12} \\
&\quad{}+ 83317\,x^2 y^{13} - 209\,x y^{14} + 1443\,x^{14} y^2 - 316732\,x^{13} y^3 + 5171738\,x^{12} y^4 \\
&\quad{}+ 39864396\,x^{11} y^5 + 60304101\,x^{10} y^6 - 24538800\,x^9 y^7 - 99382140\,x^8 y^8 \\
&\quad{}- 24538800\,x^7 y^9 + 60304101\,x^6 y^{10} + 39864396\,x^5 y^{11} + 5171738\,x^4 y^{12} \\
&\quad{}- 316732\,x^3 y^{13} + 1443\,x^2 y^{14} - 6110\,x^{14} y^3 + 814385\,x^{13} y^4 - 6605027\,x^{12} y^5 \\
&\quad{}- 55000374\,x^{11} y^6 - 122448612\,x^{10} y^7 - 146619330\,x^9 y^8 - 146619330\,x^8 y^9 \\
&\quad{}- 122448612\,x^7 y^{10} - 55000374\,x^6 y^{11} - 6605027\,x^5 y^{12} + 814385\,x^4 y^{13} - 6110\,x^3 y^{14} \\
&\quad{}+ 17732\,x^{14} y^4 - 1493206\,x^{13} y^5 + 4863144\,x^{12} y^6 + 49479144\,x^{11} y^7 \\
&\quad{}+ 116068524\,x^{10} y^8 + 146005860\,x^9 y^9 + 116068524\,x^8 y^{10} + 49479144\,x^7 y^{11} \\
&\quad{}+ 4863144\,x^6 y^{12} - 1493206\,x^5 y^{13} + 17732\,x^4 y^{14} - 37323\,x^{14} y^5 + 2003001\,x^{13} y^6 \\
&\quad{}- 646932\,x^{12} y^7 - 26836524\,x^{11} y^8 - 59751978\,x^{10} y^9 - 59751978\,x^9 y^{10} \\
&\quad{}- 26836524\,x^8 y^{11} - 646932\,x^7 y^{12} + 2003001\,x^6 y^{13} - 37323\,x^5 y^{14} + 58773\,x^{14} y^6 \\
&\quad{}- 1980264\,x^{13} y^7 - 2799654\,x^{12} y^8 + 5599308\,x^{11} y^9 + 12392094\,x^{10} y^{10} \\
&\quad{}+ 5599308\,x^9 y^{11} - 2799654\,x^8 y^{12} - 1980264\,x^7 y^{13} + 58773\,x^6 y^{14} - 70356\,x^{14} y^7 \\
&\quad{}+ 1428999\,x^{13} y^8 + 3366649\,x^{12} y^9 + 2877732\,x^{11} y^{10} + 2877732\,x^{10} y^{11} \\
&\quad{}+ 3366649\,x^9 y^{12} + 1428999\,x^8 y^{13} - 70356\,x^7 y^{14} + 64350\,x^{14} y^8 - 727870\,x^{13} y^9 \\
&\quad{}- 2000284\,x^{12} y^{10} - 2455908\,x^{11} y^{11} - 2000284\,x^{10} y^{12} - 727870\,x^9 y^{13} + 64350\,x^8 y^{14} \\
&\quad{}- 44759\,x^{14} y^9 + 240383\,x^{13} y^{10} + 662818\,x^{12} y^{11} + 662818\,x^{11} y^{12} + 240383\,x^{10} y^{13} \\
&\quad{}- 44759\,x^9 y^{14} + 23309\,x^{14} y^{10} - 38428\,x^{13} y^{11} - 99034\,x^{12} y^{12} - 38428\,x^{11} y^{13} \\
&\quad{}+ 23309\,x^{10} y^{14} - 8814\,x^{14} y^{11} - 3653\,x^{13} y^{12} - 3653\,x^{12} y^{13} - 8814\,x^{11} y^{14} \\
&\quad{}+ 2288\,x^{14} y^{12} + 2758\,x^{13} y^{13} + 2288\,x^{12} y^{14} - 365\,x^{14} y^{13} - 365\,x^{13} y^{14} + 27\,x^{14} y^{14})\\
&\;{}/[(2 - x - y + x y)^{27} (4\,x + 4\,y - 4\,x^2 - 10\,x y - 4\,y^2 + x^3 + 10\,x^2 y + 10\,x y^2 + y^3 \\
&\quad{}- 4\,x^3 y - 11\,x^2 y^2 - 4\,x y^3 + 5\,x^3 y^2 + 5\,x^2 y^3 - 2\,x^3 y^3)].
\end{align*}
\end{small}
\!\!Notice that the symmetry of $f_{15}$ can be verified by inspection.

The formula for $f_{21}^*$ in Theorem~\ref{Thm:windiff-mr} is too complicated to include here, but its restriction to the diagonal $y=x$ is given by
\begin{small}
\begin{align*}
&f_{21}^*(x,x) \\
&{}=-(1-2\,x)(4 - 118\,x + 2965\,x^2 - 57152\,x^3 + 868490\,x^4 - 10716456\,x^5 \\
&\quad{}+ 109968960\,x^6 - 956116176\,x^7 + 7146363372\,x^8 - 46448967196\,x^9 \\
&\quad{}+ 264950105372\,x^{10} - 1336163457956\,x^{11} + 5993324105124\,x^{12} \\
&\quad{}- 24027439349760\,x^{13} + 86436662246164\,x^{14} - 279912771953512\,x^{15} \\
&\quad{}+ 818043980666416\,x^{16} - 2161728799462794\,x^{17} + 5172607725929342\,x^{18} \\
&\quad{}- 11217711793111134\,x^{19} + 22059150913752150\,x^{20} - 39334435754195620\,x^{21} \\
&\quad{}+ 63573102600730426\,x^{22} - 93048429757800920\,x^{23} + 123165852823525376\,x^{24} \\
&\quad{}- 147164328676371788\,x^{25} + 158336669226676692\,x^{26} - 152925254595431252\,x^{27} \\
&\quad{}+ 132073949699924276\,x^{28} - 101513103939969032\,x^{29} + 69029907659206580\,x^{30} \\
&\quad{}- 41228526090855584\,x^{31} + 21431433409681832\,x^{32} - 9584983884534944\,x^{33} \\
&\quad{}+ 3633685742908742\,x^{34} - 1144761243925196\,x^{35} + 291604321800992\,x^{36} \\
&\quad{}- 57699710225024\,x^{37} + 8319988540052\,x^{38} - 777595791224\,x^{39} \\
&\quad{}+ 35345268560\,x^{40} + 9520\,x^{41} + 4760\,x^{42} + 1360\,x^{43} + 2720\,x^{44} + 1360\,x^{45} \\
&\quad{}+ 272\,x^{46} + 544\,x^{47} + 272\,x^{48} + 34\,x^{49} + 68\,x^{50} + 34\,x^{51} + 2\,x^{52} + 4\,x^{53} + 2\,x^{54}) \\
&/[(2-x)(4-3\,x)].
\end{align*}
\end{small}

The formula for $f_{21}^\circ$ in Theorem~\ref{Thm:windiff-hr} is too complicated to include here, but its restriction to the diagonal $y=x$ is given by
\begin{small}
\begin{align*}
&f_{21}^\circ(x,x) \\
&{}=(16 - 840\,x + 14932\,x^2 - 60826\,x^3 - 1657304\,x^4 + 30344504\,x^5 - 234880753\,x^6 \\
&\quad{}+ 786754474\,x^7 + 1909967494\,x^8 - 35455692788\,x^9 + 189290132070\,x^{10} \\
&\quad{}- 506157712476\,x^{11} + 121451108476\,x^{12} + 4896741627146\,x^{13} \\
&\quad{}- 22369020701380\,x^{14} + 53754181367220\,x^{15}  - 63834799507644\,x^{16} \\
&\quad{}- 47334758957000\,x^{17} + 396268138740360\,x^{18} - 935229590903672\,x^{19} \\
&\quad{}+ 1252746627591250\,x^{20} - 660830414520276\,x^{21} - 1253215350015590\,x^{22} \\
&\quad{}+ 3965304099511750\,x^{23} - 5994285289123358\,x^{24} + 5863204115149160\,x^{25} \\
&\quad{}- 3360398642443608\,x^{26} - 101264695791000\,x^{27} + 2491651099450144\,x^{28} \\
&\quad{}- 2674134865260578\,x^{29} + 1075700510092416\,x^{30} + 942990102862852\,x^{31} \\
&\quad{}- 2243017295261572\,x^{32} + 2533591593132952\,x^{33} - 2173253860525280\,x^{34} \\
&\quad{}+ 1633316631069820\,x^{35} - 1171929780703780\,x^{36} + 834995815965076\,x^{37} \\
&\quad{}- 587328832064166\,x^{38} + 398414099101024\,x^{39} - 256301511885384\,x^{40} \\
&\quad{}+ 155326261176392\,x^{41} - 88521994339764\,x^{42} + 47423943506656\,x^{43} \\
&\quad{}- 23881795684112\,x^{44} + 11295318189792\,x^{45} - 5002059387432\,x^{46} \\
&\quad{}+ 2063758905696\,x^{47} - 789951928976\,x^{48} + 280085379928\,x^{49} - 91896301464\,x^{50} \\
&\quad{}+ 27751064996\,x^{51} - 7618506990\,x^{52} + 1872431576\,x^{53}- 407928788\,x^{54} \\
&\quad{}+ 78960018\,x^{55} - 13744626\,x^{56} + 2160136\,x^{57} - 297916\,x^{58} + 33704\,x^{59} \\
&\quad{}- 2836\,x^{60} + 154\,x^{61} - 4\,x^{62}) \\
&{}/[(2 -x) (2 - 2\,x +\,x^2) (8 - 18\,x + 22\,x^2 - 19\,x^3 + 10\,x^4 - 2\,x^5)].
\end{align*}
\end{small}

\end{document}